\documentclass[12pt]{article}
\usepackage{rawfonts,graphicx,amssymb,amsmath,amsfonts,t1enc}
\input{prepictex}
\input{pictex}
\input{postpictex}
\newcommand{\R}{\mathbb R}
\newcommand{\T}{\textbf{\textit{T}}}

\begin{document}
\hfill
\centerline{Steinhaus' problem on partition of a triangle}
\vskip 0.1truecm
\par
\noindent
\centerline{\footnotesize (to appear in Forum Geometricorum)}
\vskip 0.1truecm
\par
\noindent
\centerline{Apoloniusz Tyszka}
\vskip 0.1truecm
\par
\noindent
{\bf Abstract.} H. Steinhaus asked a question whether inside each
acute triangle there is a point from which perpendiculars to the sides
divide the triangle into three parts with equal areas.
We present two methods of solving Steinhaus' problem.
\vskip 0.1truecm
\par
The $n$-dimensional case of Theorem 1 was proved in \cite{Kuratowski},
see also \cite{Borsuk} and [3, Theorem~2.1, p.~152]. For an earlier mass-partition
version of Theorem~1, for bounded convex masses in ${\R}^n$ and
$r_1=r_2=...=r_{n+1}$, see \cite{Levi}.
\vskip 0.1truecm
\par
\noindent
{\bf Theorem 1.} Let $\T \subseteq {\R}^2$ be a bounded measurable set,
and let $|\T|$ be the measure of~$\T$.
Let $\alpha_1$, $\alpha_2$, $\alpha_3$ be the angles determined
by three rays emanating from a point, and let $\measuredangle \alpha_1<\pi$,
$\measuredangle \alpha_2<\pi$, $\measuredangle \alpha_3<\pi$. Let $r_1$, $r_2$, $r_3$ be
non-negative numbers such that $r_1+r_2+r_3=|\T|$.
Then there exists a translation
$\lambda: {\R}^2 \to {\R}^2$ such that $|\lambda(\T) \cap \alpha_1|=r_1$,
$|\lambda(\T)\cap \alpha_2|=r_2$, $|\lambda(\T)\cap \alpha_3|=r_3$.
\vskip 0.1truecm
\par
H. Steinhaus asked a question (\cite{Steinhaus1},~\cite{Steinhaus2})
whether inside each acute triangle there is a point from which
perpendiculars to the sides divide the triangle into three parts
with equal areas. Theorem~1 indicates only that there is such
a point in the plane. There is no proof that partitions such as
presented in Figure~1 and Figure~2 are not partitions into three
parts with equal areas.
\vskip 0.1truecm
\par
\noindent
\centerline{
\beginpicture
\setcoordinatesystem units <0.1 cm, 0.1 cm>
\setplotarea x from 0 to 125, y from -20 to 25
\setplotsymbol({ .})
\plot 0 0 45 0 /
\plot 45 0 20 25 /
\plot 0 0 20 25 /
\plot 80 0 125 0 /
\plot 125 0 100 25 /
\plot 80 0 100 25 /
\plot 25 -15 25 -5 /
\plot 25 -5 4.19 11.65 /
\plot 25 -5 40.33 10.33 /
\plot 105 -15 105 0 /
\plot 105 0 86.63 14.69 /
\plot 105 0 117.83 12.83 /
\circulararc -90 degrees from 9.82 12.27 center at 7.32 9.15
\circulararc 90 degrees from 34.67 10.33 center at 37.5 7.5
\circulararc -90 degrees from 25 -10 center at 25 0
\circulararc -90 degrees from 92.25 15.32 center at 89.76 12.2
\circulararc 90 degrees from 112.17 12.83 center at 115 10
\circulararc -90 degrees from 105 -4 center at 105 0
\put {Figure 1} at 23 -20
\put {Figure 2} at 103 -20
\put {\large \bf .} at 10.25 9.5
\put {\large \bf .} at 35.75 7.5
\put {\large \bf .} at 20 -5
\put {\large \bf .} at 104 -1.75
\put {\large \bf .} at 92.75 12.25
\put {\large \bf .} at 113.25 10
\endpicture}
\vskip 0.1truecm
\par
Long and elementary solutions of Steinhaus' problem appeared in:
[7, pp.~101--104], [8, pp.~103--105], [11, pp.~133--138], \cite{Stojda}.
For some acute triangle with rational coordinates of vertices,
the point solving Steinhaus' problem is not constructible with
ruler and compass alone, see \cite{Tyszka2}.
Following article \cite{Tyszka1},
we will present two methods of solving Steinhaus' problem.
An anonymous referee found another two solutions which are constructive.
His article will appear in this journal.
\vskip 0.2truecm
\par
For $X \in \triangle ABC$ by $P(A,X)$, $P(B,X)$, $P(C,X)$
we denote figure areas at the vertices $A$, $B$, $C$ respectively,
see Figure~3.
The areas $P(A,X)$, $P(B,X)$, $P(C,X)$ are continuous as functions
of $X$ on the triangle $ABC$. The function
$f(X)={\rm min}\{P(A,X),P(B,X),P(C,X)\}$ is also continuous.
By Weierstrass' theorem $f$ attains a supremum on the triangle~$ABC$
i.e. there exists $X_0 \in \triangle ABC$ such that $f(X) \leq f(X_0)$
for all $X \in \triangle ABC$.
\vskip 0.2truecm
\par
\noindent
\centerline{
\beginpicture
\setcoordinatesystem units <0.1 cm, 0.1 cm>
\setplotarea x from -52 to 0, y from -3 to 46
\setplotsymbol({ .})
\plot -50 14 0 0 /
\plot 0 0 0 43 /
\plot 0 43 -50 14 /
\plot -20 20 -23.74 6.65 /
\plot -20 20 0 20 /
\plot -20 20 -24.95 28.53 /
\circulararc -90 degrees from -27.59 7.73 center at -23.74 6.65
\circulararc 90 degrees from -4 20 center at 0 20
\circulararc 90 degrees from -22.94 25.07 center at -24.95 28.53
\put {$A$} at -52 14
\put {$B$} at 0 -3
\put {$C$} at 0 46
\put {$X$} at -17 17
\put {$P(A,X)$} at -34 15
\put {$P(B,X)$} at -8 7
\put {$P(C,X)$} at -8 30
\put {\large \bf .} at -24.2 8.7
\put {\large \bf .} at -1 18.5
\put {\large \bf .} at -22 28
\endpicture}
\vskip 0.2truecm
\par
\noindent
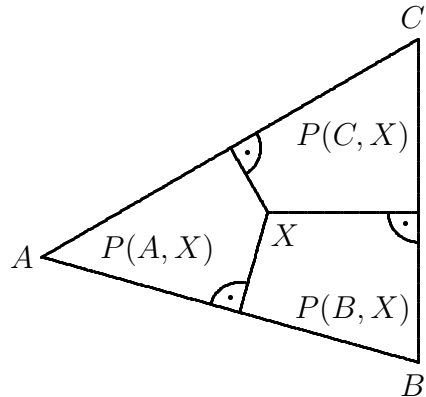
\centerline{Figure 3}
\vskip 0.2truecm
\par
\noindent
{\bf Lemma.} For a point $X$ lying on a side of an acute triangle,
the area at the opposite vertex is greater than some other area.
\vskip 0.2truecm
\par
\noindent
{\it Proof.} Without loss of generality, we may assume that
$X \in \overline{AB}$ and $|AX| \leq |BX|$, see Figure~4.
Straight line~$XX'$ parallel to straight line~$BC$ cuts the
triangle~$AXX'$ greater than~$P(A,X)$ (as the angle~$ACB$ is acute),
but not greater than the triangle $CXX'$ because
$|AX'|<\frac{|AC|}{2}<|X'C|$. Hence
$P(A,X)<|\triangle AXX'| \leq |\triangle CXX'|<P(C,X)$.
\vskip 0.2truecm
\par
\noindent
\centerline{
\beginpicture
\setcoordinatesystem units <0.1 cm, 0.1 cm>
\setplotarea x from -2 to 85, y from -3 to 48
\setplotsymbol({ .})
\plot 0 0 82 0 /
\plot 82 0 37 45 /
\plot 37 45 0 0 /
\plot 33 0 57.5 24.5 /
\plot 33 0 13.31 16.19 /
\plot 33 0 34.51 16.93 /
\plot 37 45 35.05 23.09 /
\circulararc 90 degrees from 10.77 13.1 center at 13.31 16.19
\circulararc 90 degrees from 54.67 21.67 center at 57.5 24.5
\setdashes
\plot 33 0 14.89 18.11 /
\put {$A$} at -2 0
\put {$B$} at 85 0
\put {$C$} at 39 48
\put {$X$} at 33 -3
\put {$P(A,X)$} at 16 5
\put {$P(B,X)$} at 57 5
\put {$P(C,X)$} at 35 19.5
\put {\large \bf .} at 58 22
\put {\large \bf .} at 14 14
\endpicture}
\vskip 0.2truecm
\par
\noindent
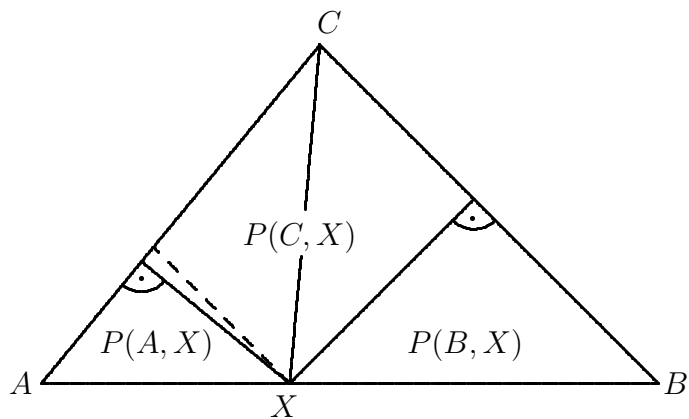
\centerline{Figure 4}
\newline
\rightline{$\Box$}
\newpage
\par
\noindent
{\bf Theorem 2.} If a triangle $ABC$ is acute and $f$ attains a supremum at $X_0$,
then $P(A,X_0)=P(B,X_0)=P(C,X_0)=\frac{|\triangle ABC|}{3}$.
\vskip 0.2truecm
\par
\noindent
{\it Proof.} $f(A)=f(B)=f(C)=0$, and $0$ is not a supremum of~$f$.
Therefore $X_0$ is not a vertex of the triangle~$ABC$.
Let us assume, to set the attention, that $f(X_0)=P(A,X_0)$.
By the Lemma $X_0 \not\in \overline{BC}$.
Suppose, on the contrary, that some of the other areas,
let's say $P(C,X_0)$, is greater than $P(A,X_0)$.
\par
Case 1: $X_0 \not\in \overline{AC}$.
When shifting $X_0$ from the segment~$\overline{AB}$ by appropriately
small~$\varepsilon$ and perpendicularly to the segment~$\overline{AB}$
(see Figure~5), we receive $P(C,X)$ further greater than $f(X_0)$
and at the same time $P(A,X)>P(A,X_0)$ and $P(B,X)>P(B,X_0)$.
Hence $f(X)>f(X_0)$, a contradiction.
\par
\noindent
\centerline{
\beginpicture
\setcoordinatesystem units <0.08 cm, 0.08 cm>
\setplotarea x from -2 to 61, y from -23 to 38
\setplotsymbol({ .})
\plot 0 0 61 -20 /
\plot 61 -20 61 35 /
\plot 0 0 61 35 /
\plot 37 3 32.52 -10.66 /
\plot 37 3 29.13 16.71 /
\plot 37 3 61 3 /
\circulararc 90 degrees from 25.66 14.72 center at 29.13 16.71
\circulararc 90 degrees from 36.32 -11.91 center at 32.52 -10.66
\circulararc 90 degrees from 57 3 center at 61 3
\arrow <2mm> [0.6, 1.5] from 37 3 to 38.87 8.7
\setdashes
\plot 38.87 8.7 33 18.93 /
\plot 38.87 8.7 61 8.7 /
\put {$A$} at -2 0
\put {$B$} at 61 -23
\put {$C$} at 61 38
\put {$X_0$} at 40.5 -0.5
\put {$X$} at 42.5 12
\put {$\varepsilon$} at 43 6
\put {\large \bf .} at 35.5 -9.5
\put {\large \bf .} at 29.5 14.5
\put {\large \bf .} at 60.25 1.25
\endpicture}
\par
\noindent
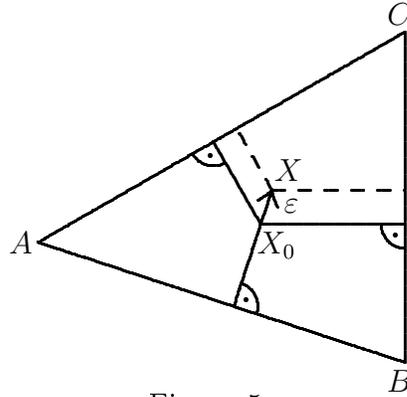
\centerline{Figure 5}
\vskip 0.2truecm
\par
Case 2: $X_0 \in \overline{AC} \setminus \{A,C\}$. By the Lemma
\par
\noindent
\centerline{$P(B,X_0)>{\rm min}\{P(A,X_0),P(C,X_0)\} \geq {\rm min} \{P(A,X_0),P(B,X_0),P(C,X_0)\}=f(X_0)$}
\par
\noindent
When shifting $X_0$ from the segment~$\overline{AC}$ by appropriately
small~$\varepsilon$ and perpendicularly to the segment~$\overline{AC}$
(see Figure~6), we receive $P(B,X)$ further greater than $f(X_0)$
and at the same time $P(A,X)>P(A,X_0)$ and $P(C,X)>P(C,X_0)$.
Hence $f(X)>f(X_0)$, a contradiction.
\par
\noindent
\centerline{
\beginpicture
\setcoordinatesystem units <0.08 cm, 0.08 cm>
\setplotarea x from -2 to 61, y from -23 to 38
\setplotsymbol({ .})
\plot 0 0 61 -20 /
\plot 61 -20 61 35 /
\plot 0 0 61 35 /
\plot 26.02 14.93 61 14.93 /
\plot 26.02 14.93 19.07 -6.25 /
\arrow <2mm> [0.6, 1.5] from 26.02 14.93 to 30 8
\circulararc 90 degrees from 61 18.93 center at 61 14.93
\circulararc 90 degrees from 20.32 -2.45 center at 19.07 -6.25
\circulararc -90 degrees from 30 8 center at 26.02 14.93
\setdashes
\plot 30 8 61 8 /
\plot 30 8 24.72 -8.1 /
\put {$A$} at -2 0
\put {$B$} at 61 -23
\put {$C$} at 61 38
\put {$X_0$} at 24 17.75
\put {$X$} at 33 5
\put {$\varepsilon$} at 30.5 13
\put {\large \bf .} at 18.75 -4
\put {\large \bf .} at 60.25 16.75
\put {\large \bf .} at 23 10
\endpicture}
\par
\noindent
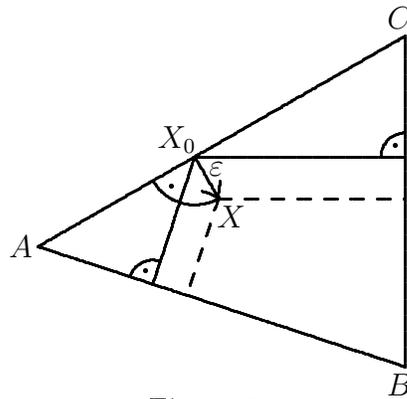
\centerline{Figure 6}
\newline
\rightline{$\Box$}
\vskip 0.2truecm
\par
For each acute triangle $ABC$ there is a unique
$X_0 \in \triangle ABC$ such that
$P(A,X_0)=P(B,X_0)=P(C,X_0)=\frac{|\triangle ABC|}{3}$.
Indeed, if $X \neq X_0$ then $X$ lies in some of the
quadrangles determined by~$X_0$. Let's say that $X$ lies
in the quadrangle with vertex~$A$, see Figure~7.
Then $P(A,X)<P(A,X_0)=\frac{|\triangle ABC|}{3}$.
\vskip 0.2truecm
\par
\noindent
\centerline{
\beginpicture
\setcoordinatesystem units <0.1 cm, 0.1 cm>
\setplotarea x from -52 to 0, y from -3 to 46
\setplotsymbol({ .})
\plot -50 14 0 0 /
\plot 0 0 0 43 /
\plot 0 43 -50 14 /
\plot -20 20 -23.74 6.65 /
\plot -20 20 0 20 /
\plot -20 20 -24.95 28.53 /
\circulararc 90 degrees from -19.89 5.57 center at -23.74 6.65
\circulararc 90 degrees from -4 20 center at 0 20
\circulararc 90 degrees from -22.94 25.07 center at -24.95 28.53
\setdashes
\plot -25 23 -29.15 8.16 /
\plot -25 23 -27.39 27.12 /
\plot -25 23 0 23 /
\put {$A$} at -52 14
\put {$B$} at 0 -3
\put {$C$} at 0 46
\put {$X_0$} at -17 17
\put {$X$} at -27.5 23
\put {\large \bf .} at -21 8
\put {\large \bf .} at -1 18.5
\put {\large \bf .} at -22 28
\endpicture}
\vskip 0.2truecm
\par
\noindent
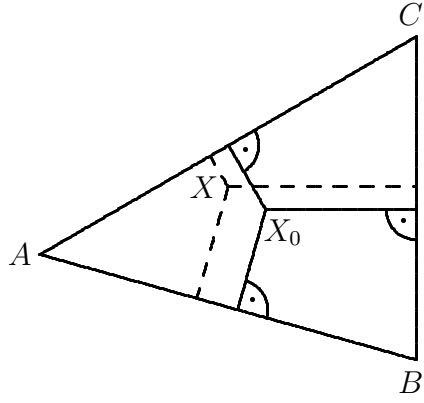
\centerline{Figure 7}
\vskip 0.2truecm
\par
The sets $R_A=\{X \in \triangle ABC: P(A,X)=f(X)\}$,
$R_B=\{X \in \triangle ABC: P(B,X)=f(X)\}$
and $R_C=\{X \in \triangle ABC: P(C,X)=f(X)\}$
are closed and cover the triangle $ABC$.
Assume that the triangle $ABC$ is acute.
By the Lemma $R_A \cap \overline{BC}=\emptyset$,
$R_B \cap \overline{AC}=\emptyset$ and $R_C \cap \overline{AB}=\emptyset$.
The theorem proved in \cite{Knaster} guarantees that
$R_A \cap R_B \cap R_C \neq \emptyset$, see also [3, item~D4, p.~101]
and [1, item~2.23, p.~162]. Any point belonging to $R_A \cap R_B \cap R_C$
lies inside the triangle~$ABC$ and determines the partition of the
triangle~$ABC$ into three parts with equal areas.
\vskip 0.2truecm
\par
The above proof remains valid for all right triangles, because
the thesis of the Lemma holds for all right triangles.
For each triangle we have:
\begin{description}
\item{{\bf 1)}}
There is a unique point in the plane which determines
the partition of the triangle into three equal areas.
\item{{\bf 2)}}
The point of partition into three equal areas lies inside
the triangle if and only if the thesis of the Lemma holds for the triangle.
\item{{\bf 3)}}
The point of partition into three equal areas lies inside
the triangle if and only if the supremum of $f$ on the boundary
of the triangle is smaller than the supremum of $f$ on the whole triangle.
For each acute or right triangle $ABC$, the supremum of $f$ on
the boundary does not exceed $\frac{|\triangle ABC|}{4}$.
\newpage
\item{{\bf 4)}}
The point of partition into three equal areas lies inside the
triangle, if the triangle has two angles in the interval
$({\rm arctg}(\frac{1}{\sqrt{2}}), \frac{\pi}{2}]$.
This condition holds for each acute or right triangle.
\item{{\bf 5)}}
If the point of partition into three equal areas lies
inside the triangle, then it is a partition into quadrangles.
\end{description}
\vskip 0.2truecm
\par
For each triangle $ABC$ with $\measuredangle C > \frac{\pi}{2}$ the point
of partition into three equal areas lies inside the triangle if and
only if
$$
\sqrt{(1+{\rm {tg}^2~} \measuredangle A) \cdot {\rm tg~} \measuredangle B}+
\sqrt{(1+{\rm {tg}^2~} \measuredangle B) \cdot {\rm tg~} \measuredangle A}
>\sqrt{3 \cdot ({\rm tg~} \measuredangle A + {\rm tg~} \measuredangle B)}
$$
In case $\measuredangle C > \frac{\pi}{2}$ and
$$
\sqrt{(1+{\rm {tg}^2~} \measuredangle A) \cdot {\rm tg~} \measuredangle B}+
\sqrt{(1+{\rm {tg}^2~} \measuredangle B) \cdot {\rm tg~} \measuredangle A}
=\sqrt{3 \cdot ({\rm tg~} \measuredangle A + {\rm tg~} \measuredangle B)}
$$
the unique $X_0 \in \overline{AB}$ such that
$|AX_0|=\sqrt{\frac{(1+{\rm {tg}^2} \measuredangle A) \cdot {\rm tg~} \measuredangle B}
{3 \cdot ({\rm tg~} \measuredangle A+{\rm tg~} \measuredangle B)}} \cdot |AB|$
and
$|BX_0|=\sqrt{\frac{(1+{\rm {tg}^2} \measuredangle B) \cdot {\rm tg~} \measuredangle A}
{3 \cdot ({\rm tg~} \measuredangle A+{\rm tg~} \measuredangle B)}} \cdot |AB|$
determines the partition of the triangle $ABC$ into three equal areas.
It is a partition into a triangle with vertex $A$, and a triangle with vertex $B$,
and a quadrangle. There remains the case when $\measuredangle C > \frac{\pi}{2}$ and
\begin{equation}
\tag*{$(\ast)$}
\sqrt{(1+{\rm {tg}^2~} \measuredangle A) \cdot {\rm tg~} \measuredangle B}+
\sqrt{(1+{\rm {tg}^2~} \measuredangle B) \cdot {\rm tg~} \measuredangle A}
<\sqrt{3 \cdot ({\rm tg~} \measuredangle A + {\rm tg~} \measuredangle B)}
\end{equation}
There is a straight line $a$ perpendicular to the segment
$\overline{AC}$ which cuts from the triangle $ABC$ a figure with
the area $\frac{|\triangle ABC|}{3}$, see Figure~8. There is a straight line $b$
perpendicular to the segment $\overline{BC}$ which cuts from
the triangle $ABC$ a figure with the area $\frac{|\triangle ABC|}{3}$,
see Figure~8. By $(\ast)$, the intersection point of straight lines~$a$
and~$b$ lies outside the triangle $ABC$. This point determines the
partition of the triangle $ABC$ into three equal areas.
\vskip 0.2truecm
\par
\noindent
\centerline{
\beginpicture
\setcoordinatesystem units <0.1 cm, 0.1 cm>
\setplotarea x from -18 to 37.5, y from -15.5 to 23
\setplotsymbol({ .})
\plot 0 0 37 0 /
\plot 0 0 -13 13 /
\plot 37 0 -13 13 /
\plot 12.17 -15 12.17 14.8 /
\plot -17.63 -15 12.17 14.8 /
\plot 12.17 14.8 14.2 22.62 /
\put {.} at 5 7
\put {.} at 5 8
\put {.} at 4 6
\put {.} at 4 7
\put {.} at 4 8
\put {.} at 3 5
\put {.} at 3 6
\put {.} at 3 7
\put {.} at 3 8
\put {.} at 3 9
\put {.} at 2 4
\put {.} at 2 5
\put {.} at 2 6
\put {.} at 2 7
\put {.} at 2 8
\put {.} at 2 9
\put {.} at 1 3
\put {.} at 1 4
\put {.} at 1 5
\put {.} at 1 6
\put {.} at 1 7
\put {.} at 1 8
\put {.} at 1 9
\put {.} at 0 2
\put {.} at 0 3
\put {.} at 0 4
\put {.} at 0 5
\put {.} at 0 6
\put {.} at 0 7
\put {.} at 0 8
\put {.} at 0 9
\put {.} at -1 2
\put {.} at -1 3
\put {.} at -1 4
\put {.} at -1 5
\put {.} at -1 6
\put {.} at -1 7
\put {.} at -1 8
\put {.} at -1 9
\put {.} at -1 10
\put {.} at -2 3
\put {.} at -2 4
\put {.} at -2 5
\put {.} at -2 6
\put {.} at -2 7
\put {.} at -2 8
\put {.} at -2 9
\put {.} at -2 10
\put {.} at -3 4
\put {.} at -3 5
\put {.} at -3 6
\put {.} at -3 7
\put {.} at -3 8
\put {.} at -3 9
\put {.} at -3 10
\put {.} at -4 5
\put {.} at -4 6
\put {.} at -4 7
\put {.} at -4 8
\put {.} at -4 9
\put {.} at -4 10
\put {.} at -5 6
\put {.} at -5 7
\put {.} at -5 8
\put {.} at -5 9
\put {.} at -5 10
\put {.} at -5 11
\put {.} at -6 7
\put {.} at -6 8
\put {.} at -6 9
\put {.} at -6 10
\put {.} at -6 11
\put {.} at -7 8
\put {.} at -7 9
\put {.} at -7 10
\put {.} at -7 11
\put {.} at -8 9
\put {.} at -8 10
\put {.} at -8 11
\put {.} at -8 12
\put {.} at -9 10
\put {.} at -9 11
\put {.} at -9 12
\put {.} at -10 11
\put {.} at -10 12
\put {.} at -11 12
\put {.} at 14 1
\put {.} at 14 2
\put {.} at 14 3
\put {.} at 14 4
\put {.} at 14 5
\put {.} at 14 6
\put {.} at 15 1
\put {.} at 15 2
\put {.} at 15 3
\put {.} at 15 4
\put {.} at 15 5
\put {.} at 16 1
\put {.} at 16 2
\put {.} at 16 3
\put {.} at 16 4
\put {.} at 16 5
\put {.} at 17 1
\put {.} at 17 2
\put {.} at 17 3
\put {.} at 17 4
\put {.} at 17 5
\put {.} at 18 1
\put {.} at 18 2
\put {.} at 18 3
\put {.} at 18 4
\put {.} at 18 5
\put {.} at 19 1
\put {.} at 19 2
\put {.} at 19 3
\put {.} at 19 4
\put {.} at 20 1
\put {.} at 20 2
\put {.} at 20 3
\put {.} at 20 4
\put {.} at 21 1
\put {.} at 21 2
\put {.} at 21 3
\put {.} at 21 4
\put {.} at 22 1
\put {.} at 22 2
\put {.} at 22 3
\put {.} at 22 4
\put {.} at 23 1
\put {.} at 23 2
\put {.} at 23 3
\put {.} at 24 1
\put {.} at 24 2
\put {.} at 24 3
\put {.} at 25 1
\put {.} at 25 2
\put {.} at 25 3
\put {.} at 26 1
\put {.} at 26 2
\put {.} at 27 1
\put {.} at 27 2
\put {.} at 28 1
\put {.} at 28 2
\put {.} at 29 1
\put {.} at 29 2
\put {.} at 30 1
\put {.} at 30 2
\put {.} at 31 1
\put {.} at 32 1
\put {.} at 33 1
\circulararc 90 degrees from 12.17 -4 center at 12.17 0
\circulararc 90 degrees from -4.15 4.15 center at -1.32 1.32
\circulararc -90 degrees from 12.17 14.8 center at 10.14 6.98
\setdashes
\plot 12.17 14.8 10.14 6.98 /
\put {\large \bf .} at 14.5 -1.5
\put {\large \bf .} at -3 1.5
\put {\large \bf .} at 15 9.5
\put {$b$} at 15 -10
\put {$a$} at -15 -10
\put {$A$} at -14.5 13
\put {$C$} at 0 -3
\put {$B$} at 37 -3
\endpicture}
\vskip 0.2truecm
\par
\noindent
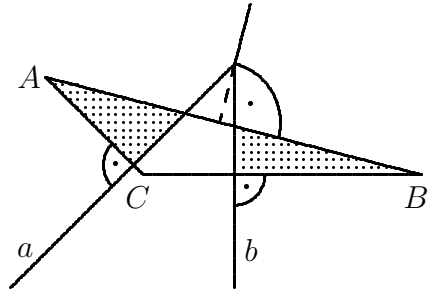
\centerline{Figure 8}

Apoloniusz Tyszka\\
Technical Faculty\\
Hugo Ko\l{}\l{}\k{a}taj University\\
Balicka 116B, 30-149 Krak\'ow, Poland\\
E-mail address: {\it rttyszka@cyf-kr.edu.pl}

\begin{thebibliography}{14}
\bibitem{Alexandrov}
{\sc P. S. Alexandrov}, {\it Combinatorial topology},
Dover Publications, Mineola, NY, 1998.
\bibitem{Borsuk}
{\sc K. Borsuk}, {\it An application of the theorem on antipodes to
the measure theory}, Bull. Acad. Polon. Sci. Cl. III 1~(1953), pp.~87--90.
\bibitem{Granas}
{\sc A. Granas and J. Dugundji}, {\it Fixed point theory},
Springer-Verlag, New York, 2003.
\bibitem{Knaster}
{\sc B. Knaster, K. Kuratowski and S. Mazurkiewicz},
{\it Ein Beweis des Fixpunktsatzes fur $n$-dimensionale Simplexe},
Fund. Math. 14~(1929), pp.~132--137;
reprinted in: K. Kuratowski, {\it Selected papers},
PWN (Polish Scientific Publishers), Warsaw, 1988, pp.~332--337,
S. Mazurkiewicz, {\it Travaux de topologie et ses applications},
PWN (\'Editions Scientifiques de Pologne), Warsaw, 1969, pp.~192--197.
\bibitem{Kuratowski}
{\sc K. Kuratowski and H. Steinhaus},
{\it Une application g\'eom\'etrique du th\'eor\'eme de Brouwer sur les
points invariants}, Bull. Acad. Polon. Sci. Cl. III 1~(1953), pp.~83--86;
reprinted in: K. Kuratowski, {\it Selected papers},
PWN (Polish Scientific Publishers), Warsaw, 1988, pp.~520--523,
H. Steinhaus, {\it Selected papers}, PWN (Polish Scientific Publishers),
Warsaw, 1985, pp.~636--639.
\bibitem{Levi}
{\sc F. Levi}, {\it Die Drittelungskurve}, Math. Z. 31 (1930),
no.~1, pp. 339--345.
\bibitem{Piegat1}
{\sc E. Piegat}, {\it Yet 105 problems of Hugo Steinhaus} (in Polish),
Oficyna Wydawnicza GiS, Wroc\l{}aw, 2000.
\bibitem{Piegat2}
{\sc E. Piegat}, {\it Known and unknown problems of Hugo Steinhaus}
(in Polish), Oficyna Wydawnicza GiS, Wroc\l{}aw, 2005.
\bibitem{Steinhaus1}
{\sc H. Steinhaus}, {\it Problem No.~132} (in Polish), Roczniki Polskiego
Towarzystwa Matematycznego (Annales Societatis Mathematicae Polonae),
Seria II, Wiadomo\'sci Matematyczne 9~(1966), no.~1, p.~99.
\bibitem{Steinhaus2}
{\sc H. Steinhaus}, {\it Problem No.~779} (in Polish), Matematyka 19 (1966), no.~2, p.~92.
\bibitem{Steinhaus3}
{\sc H. Steinhaus}, {\it Problems and reflections} (in Russian), Mir, Moscow, 1974.
\bibitem{Stojda}
{\sc W. Stojda}, {\it A solution of Problem No.~779} (in Polish), Matematyka 21 (1968), no.~5--6, pp.~267--273.
\bibitem{Tyszka1}
{\sc A. Tyszka}, {\it A solution of a problem of Steinhaus} (in Polish),
Matematyka 49~(1996), no.~1, pp.~3--5.
\bibitem{Tyszka2}
{\sc A. Tyszka}, {\it Steinhaus' problem cannot be solved with ruler
and compass alone} (in Polish), Matematyka 49~(1996), no.~4, pp.~238--240.
\end{thebibliography}
\end{document}